\pdfoutput=1
\newif\ifpersonal 
\personaltrue % comment to remove personal notes
\RequirePackage[l2tabu,orthodox]{nag} %detect whether obsolete packages are used
\documentclass[10pt,a4paper,reqno]{amsart} %reqno places equation numbers on the right
\usepackage[normalem]{ulem}
\usepackage{amsmath,amsthm,amssymb,mathrsfs,mathtools,bm,tensor} % math related
\usepackage{ stmaryrd }
\usepackage{microtype,lmodern} % latex technical issues
\usepackage[utf8]{inputenc} % input encoding
\usepackage[T1]{fontenc} % font encoding
\usepackage{enumitem,comment,braket,xspace,tikz-cd,csquotes} % utilities
\usepackage[all,cmtip]{xy} % because the tikzcd options [shift left], [shift right] do not work on arXiv, we switched some diagrams to xymatrix
\usepackage[centering,vscale=0.7,hscale=0.8]{geometry}
\usepackage[hidelinks]{hyperref}
\usepackage[capitalize]{cleveref}
\usepackage{graphicx,scalerel}

\usepackage[utf8]{inputenc}
\usepackage{tikz-cd}
\usetikzlibrary{cd,decorations.pathreplacing}

\usepackage[mathscr]{eucal}

\theoremstyle{plain}
\newtheorem{thm-intro}{Theorem}
\newtheorem{thm}{Theorem}[section]
\newtheorem*{thm*}{Theorem}
\newtheorem{claim}[thm]{Claim}
\newtheorem{lem}[thm]{Lemma}

\newtheorem{prop}[thm]{Proposition}

\newtheorem{cor}[thm]{Corollary}

\theoremstyle{definition}
\newtheorem{defin}[thm]{Definition}

\theoremstyle{remark}
\newtheorem{rem}[thm]{Remark}
\numberwithin{equation}{section}

\usepackage{comment}
\ifpersonal
\newcommand*{\personal}[1]{\textcolor[rgb]{0.6,0.6,1}{(Personal: #1)}}

\newenvironment{newversion}{\color{purple}}{}
\usepackage[textsize=tiny]{todonotes}
\else
\newcommand*{\personal}[1]{\ignorespaces}
\excludecomment{oldversion}

\usepackage[disable]{todonotes}
\fi

\newcommand{\cC}{\mathcal C}

\DeclareFontFamily{U}{BOONDOX-calo}{\skewchar\font=45 }
\DeclareFontShape{U}{BOONDOX-calo}{m}{n}{<-> s*[1.05] BOONDOX-r-calo}{}
\DeclareFontShape{U}{BOONDOX-calo}{b}{n}{<-> s*[1.05] BOONDOX-b-calo}{}
\DeclareMathAlphabet{\mathcalboondox}{U}{BOONDOX-calo}{m}{n}
\makeatletter
\let\save@mathaccent\mathaccent
\newcommand*\if@single[3]{%
	\setbox0\hbox{${\mathaccent"0362{#1}}^H$}%
	\setbox2\hbox{${\mathaccent"0362{\kern0pt#1}}^H$}%
	\ifdim\ht0=\ht2 #3\else #2\fi
}
\newcommand*\rel@kern[1]{\kern#1\dimexpr\macc@kerna}
\newcommand*\widebar[1]{\@ifnextchar^{{\wide@bar{#1}{0}}}{\wide@bar{#1}{1}}}
\newcommand*\wide@bar[2]{\if@single{#1}{\wide@bar@{#1}{#2}{1}}{\wide@bar@{#1}{#2}{2}}}
\newcommand*\wide@bar@[3]{%
	\begingroup
	\def\mathaccent##1##2{%
		\let\mathaccent\save@mathaccent
		\if#32 \let\macc@nucleus\first@char \fi
		\setbox\z@\hbox{$\macc@style{\macc@nucleus}_{}$}%
		\setbox\tw@\hbox{$\macc@style{\macc@nucleus}{}_{}$}%
		\dimen@\wd\tw@
		\advance\dimen@-\wd\z@
		\divide\dimen@ 3
		\@tempdima\wd\tw@
		\advance\@tempdima-\scriptspace
		\divide\@tempdima 10
		\advance\dimen@-\@tempdima
		\ifdim\dimen@>\z@ \dimen@0pt\fi
		\rel@kern{0.6}\kern-\dimen@
		\if#31
		\overline{\rel@kern{-0.6}\kern\dimen@\macc@nucleus\rel@kern{0.4}\kern\dimen@}%
		\advance\dimen@0.4\dimexpr\macc@kerna
		\let\final@kern#2%
		\ifdim\dimen@<\z@ \let\final@kern1\fi
		\if\final@kern1 \kern-\dimen@\fi
		\else
		\overline{\rel@kern{-0.6}\kern\dimen@#1}%
		\fi
	}%
	\macc@depth\@ne
	\let\math@bgroup\@empty \let\math@egroup\macc@set@skewchar
	\mathsurround\z@ \frozen@everymath{\mathgroup\macc@group\relax}%
	\macc@set@skewchar\relax
	\let\mathaccentV\macc@nested@a
	\if#31
	\macc@nested@a\relax111{#1}%
	\else
	\def\gobble@till@marker##1\endmarker{}%
	\futurelet\first@char\gobble@till@marker#1\endmarker
	\ifcat\noexpand\first@char A\else
	\def\first@char{}%
	\fi
	\macc@nested@a\relax111{\first@char}%
	\fi
	\endgroup
}
\makeatother

\newcommand{\Pro}{\mathrm{Pro}}
\newcommand{\Ind}{\mathrm{Ind}}

\newcommand{\Perf}{\mathrm{Perf}}

\newcommand{\id}{\mathrm{id}}

\newcommand*{\longhookrightarrow}{\ensuremath{\lhook\joinrel\relbar\joinrel\rightarrow}}

\usetikzlibrary{decorations.markings} %arrows for open immersions and closed immersions
\tikzset{
  closed/.style = {decoration = {markings, mark = at position 0.5 with { \node[transform shape, xscale = .8, yscale=.4] {/}; } }, postaction = {decorate} },
  open/.style = {decoration = {markings, mark = at position 0.5 with { \node[transform shape, scale = .7] {$\circ$}; } }, postaction = {decorate} }
}

\DeclareMathOperator{\Fun}{Fun}

\DeclareMathOperator{\Hom}{Hom}

\DeclareMathOperator{\Map}{Map}

\DeclareMathOperator*{\colim}{colim}
\DeclareMathOperator*{\holim}{holim}

\let\originalleft\left
\let\originalright\right
\renewcommand{\left}{\mathopen{}\mathclose\bgroup\originalleft}
\renewcommand{\right}{\aftergroup\egroup\originalright}

\makeatletter
\def\DeclareMathBinOp{\@ifstar{\declaremathbinop@star}{\declaremathbinop@nostar}}
\def\declaremathbinop@star#1#2{\def#1{\test@subnexp@star#2}}
\def\test@subnexp@star#1{\@ifnextchar_{\isol@subnexp@star#1}{\test@exp@star#1}}
\def\test@exp@star#1{\@ifnextchar^{\isol@expnsub@star#1}{\mathbin{#1}}}
\def\isol@subnexp@star#1_#2{\@ifnextchar^{\eval@subnexp@star#1_#2}{\mathbin{\operatorname*{#1}_{#2}}}}
\def\eval@subnexp@star#1_#2^#3{\mathbin{\operatorname*{#1}_{#2}^{#3}}}
\def\isol@expnsub@star#1^#2{\@ifnextchar_{\eval@expnsub@star#1^#2}{\mathbin{\operatorname*{#1}^{#2}}}}
\def\eval@expnsub@star#1^#2_#3{\mathbin{\operatorname*{#1}_{#3}^{#2}}}

\def\declaremathbinop@nostar#1#2{\def#1{\test@subnexp@nostar{#2}}}
\def\test@subnexp@nostar#1{\@ifnextchar_{\isol@subnexp@nostar#1}{\test@exp@nostar#1}}
\def\test@exp@nostar#1{\@ifnextchar^{\isol@expnsub@nostar#1}{\mathbin{#1}}}
\def\isol@subnexp@nostar#1_#2{\@ifnextchar^{\eval@subnexp@nostar#1_{{#2}}}{\mathbin{\underset{#2}{#1}}}}
\def\eval@subnexp@nostar#1_#2^#3{\mathbin{\overset{#3}{\underset{#2}{#1}}}}
\def\isol@expnsub@nostar#1^#2{\@ifnextchar_{\eval@expnsub@nostar#1^{#2}}{\mathbin{\overset{#2}{#1}}}}
\def\eval@expnsub@nostar#1^#2_#3{\mathbin{\overset{#2}{\underset{#3}{#1}}}}

\def\declaremathbinop@toto#1#2{\def#1{\testoo@subnexp@toto{#2}}}
\def\testoo@subnexp@toto#1{\@ifnextchar_{\isoloo@subnexp@toto#1}{\testoo@exp@toto#1}}
\def\testoo@exp@toto#1{\@ifnextchar^{\isoloo@expnsub@toto#1}{\mathbin{#1}}}
\def\isoloo@subnexp@toto#1_#2{\@ifnextchar^{\eval@subnexp@toto#1_{{#2}}}{\mathbin{\underset{#2}{#1}}}}
\def\eval@subnexp@toto#1_#2^#3{\mathbin{\overset{#3}{\underset{#2}{#1}}}}
\def\isoloo@expnsub@toto#1^#2{\@ifnextchar_{\eval@expnsub@toto#1^{#2}}{\mathbin{\overset{#2}{#1}}}}
\def\eval@expnsub@toto#1^#2_#3{\mathbin{\overset{#2}{\underset{#3}{#1}}}}
\makeatother

\DeclareMathBinOp*{\newtimes}{\times}
\DeclareMathBinOp*{\newamalg}{\amalg}
\DeclareMathBinOp*{\newotimes}{\otimes}

\newlist{assumptions}{enumerate}{10}
\setlist[assumptions]{label*={\upshape{(\alph*)}}}

\crefname{assumptionsi}{assumption}{assumptions}
\Crefname{assumptionsi}{Assumption}{Assumptions}

\newcommand{\mN}{\mathbb N}

\newcommand{\cond}{\mathsf{Cond}}
\newcommand{\solid}{\mathsf{Solid}}
\newcommand{\zD}{\mathscr{D}}

\newcommand{\zC}{\mathscr{C}}

\newcommand{\ale}{{\aleph_0}}

\newcommand{\Tate}{\mathscr{T}\mathrm{ate}}
\DeclareMathOperator{\uhom}{\underline{Hom}}
\DeclareMathOperator{\im}{\mathrm{im}}
\DeclareMathOperator{\coker}{\mathrm{coker}}
\DeclareMathOperator{\End}{\mathrm{End}}

\begin{document}

\title{Tate modules as condensed modules}

%\author{Francesco IACCA}
%\address{Francesco IACCA, DIMAI, Firenze, Italy}
%\email{francesco.iacca@stud.unifi.it}

% \author{Andrea MAFFEI}
% \address{Andrea MAFFEI, Universit\`a di Pisa, Italy}
% \email{andrea.maffei@unipi.it }

\author{Valerio MELANI}
\address{Valerio MELANI, DIMAI, Firenze, Italy}
\email{valerio.melani@unifi.it}

\author{Hugo POURCELOT}
\address{Hugo POURCELOT, DIMAI, Firenze, Italy}
\email{hugo.pourcelot@gmail.com}

\author{Gabriele VEZZOSI}
\address{Gabriele VEZZOSI, DIMAI, Firenze, Italy}
\email{gabriele.vezzosi@unifi.it}
\date{\today}

\begin{abstract}
    We prove that the category of countable Tate modules over an arbitrary discrete ring embeds fully faithfully into that of condensed modules. If the base ring is of finite type, we characterize the essential image as generated by the free module of infinite countable rank under direct sums, duals and retracts. 
    In the $\infty$-categorical context, under the same assumption on the base ring, we establish a fully faithful embedding of the $\infty$-category of countable Tate objects in perfect complexes, with uniformly bounded tor-amplitude, into the derived $\infty$-category of condensed modules.
    The boundedness assumption is necessary to ensure fullness, as we prove via an explicit counterexample in the unbounded case. 
\end{abstract}

\maketitle

\tableofcontents
\section*{Introduction}
Tate vector spaces were first introduced in \cite{Lef} (under the name of \emph{locally linearly compact topological vector spaces}), and then reconsidered and developed in \cite{BFM} as a convenient generalization, in the context of topological vector spaces, of finite dimensional vector spaces for which some notion of \emph{dimension} and \emph{determinant} still make sense, and behave reasonably well. It is safe to say that the prototypical Tate $k$-vector space is the ring of Laurent series $k(\!(t)\!)$. Later  \cite{Drinfinite} extended this notion to \emph{families}
of Tate vector spaces over an affine scheme $\mathrm{Spec}\, R$, called Tate $R$-modules, and proved flat descent for the functor sending $R$ to the category of Tate $R$-modules. 

The approach of \cite{Drinfinite}, like in the original definition, still belongs to the realm of topological linear algebra. More recently, in \cite{BGWTateExact}, a purely categorical version (pioneered by \cite{BeiHow}) is defined and studied: the authors associate to any exact 1-category $\mathcal{C}$ another exact category $\mathsf{Tate}(\mathcal{C})$ of \emph{Tate-objects} in $\mathcal{C}$, %(whose exact structure is split if the one of $\mathcal{C}$ is), 
in such a way that, when $\mathcal{C} = \mathsf{Proj}^{\mathrm{fg}}_R$ is the category of finitely generated projective $R$-modules, $\mathsf{Tate}(\mathcal{C})$ is equivalent to the category of Tate $R$-modules defined in \cite{Drinfinite}\footnote{when cardinalities restrictions are imposed, see \cite[Thm. 5.26]{BGWTateExact} for a precise statement.}. In particular, the construction of Tate objects in an exact category $\mathcal{C}$ can be iterated to obtain the exact category $n\text{-}\mathsf{Tate}(\mathcal{C})$ of (iterated) $n$-Tate objects in $\mathcal{C}$. This is particularly useful because it allows a proper setting for the so-called \emph{normally oriented} (non-symmetric) \emph{tensor product} between (iterated) Tate vector spaces: the normally oriented tensor product of a $n$-Tate vector space and a $m$-Tate vector space is naturally a $(n+m)$-Tate vector space (see \cite{BGHW}). Moreover this construction generalizes  to (iterated) Tate objects in an exact category endowed with a suitable symmetric tensor product. We will come back to this point, as a motivation, later in this Introduction.

Several (derived) moduli stacks arising in formal loop spaces and representation theory for loop groups are expected to have cotangent complexes that are, when properly understood, possibly iterated Tate complexes (see, e.g. \cite{hennion:floops, Heleo}). This requires a derived and $\infty$-categorical version of iterated Tate objects that was developed in \cite{hennion:tate}. We will use this setting in the second part of our paper.\\

From the previous summary it is clear that the root of the theory of Tate $R$-modules has to be found in topological linear algebra. Recently, a new approach to topological algebra, vastly called Condensed Mathematics, has been proposed by D. Clausen and P. Scholze (\cite{CScond}). It is therefore natural to investigate the place of Tate $R$-modules (and, more generally, of Tate objects) inside Condensed Mathematics. This is what we try to do in this paper. \\

\noindent{\textbf{Outline of the paper.}}
We prove in Section \ref{sec:1-cat_embedding} that for any commutative ring $R$ there is a fully faithful exact functor
\[
    \mathsf{TC}_R\colon \mathsf{Tate}_{\aleph_0,R} \longhookrightarrow \mathsf{Cond}(\mathsf{Mod}_R)
\] 
where $\mathsf{Tate}_{\aleph_0,R}:=\mathsf{Tate}_{\aleph_0}(\mathsf{Proj}^{\mathrm{fg}}_R)$ is the 1-category of Tate $R$-modules whose size is countable\footnote{In the topological description of \cite{Drinfinite} these correspond to (linearly topologized) topological $R$-modules of the form $P\oplus Q^{\vee}$ where $P$ and $Q$ are discrete, countably generated $R$-modules, and $(-)^{\vee}$ denotes the topological linear dual.} (\cite[Def. 5.3 and Def. 5.23]{BGWTateExact}), and $\mathsf{Cond}(\mathsf{Mod}_R)$ is the abelian category of condensed $R$-modules, i.e. the category of accessible sheaves of $R$-modules on the pro-\'etale site of a point (\cite{CScond}). We give two equivalent descriptions of the functor $\mathsf{TC}_R$: one topological (closer to the setting of \cite{Drinfinite}) in Proposition \ref{ff1}, and the other category-theoretic (closer to \cite{BGWTateExact}) in Propositions \ref{TC'} and \ref{ff2}, and we verify that $\mathsf{TC}_R$ actually lands in the full subcategory $\mathsf{Solid}_R$ of \emph{solid} $R$-modules.

In Section \ref{sec:essential_image} we characterize its essential image in Theorem \ref{thm:ess_image}: it is the full subcategory of $\mathsf{Solid}_R$ generated by the free module $\bigoplus_{\mathbb{N}} R$ of infinite countable rank under sums, duals and retracts. Although fully faithfulness of the functor $\mathsf{TC}_R$ holds for any ring $R$, we prove this characterization of its essential image under the assumption that $R$ is of finite type over $\mathbb{Z}$.
% \begin{thm-intro}
%     [Theorem \ref{thm:ess_image}]
%     For $R$ a finitely generated $\mathbb{Z}$-algebra, the essential image of the functor 
%     $\mathrm{TS}_R \colon\mathsf{Tate}_{{\aleph}_0 , R} \to \mathsf{Solid}_R$
%     is the smallest full subcategory of $\mathsf{Solid}_R$ that contains the countable free module $\bigoplus_{\mathbb{N}} R$ and is stable under direct sums, retracts and duals.
% \end{thm-intro}

Section \ref{infinity} is devoted to extend the fully faithfulness result to the $\infty$-category of Tate $R$-modules, as defined in \cite{hennion:tate}. 
%This extension does not seem to be a trivial consequence of the 1-categorical result; in fact our approach to the $\infty$-categorical statement will be completely independent of the first two sections, avoiding the need for topological modules.  
Contrary to the 1-categorical embedding, but as in Section \ref{sec:essential_image}, we will need to assume the base ring $R$ to be of finite type over $\mathbb{Z}$. 
%Moreover, the $\infty$-categorical result turns out to hold only if we limit ourselves to perfect complexes (replacing finitely generated projective modules in the $\infty$-categorical setting) which are \emph{of fixed tor-amplitude}. 
The main result is achieved in Theorem \ref{thm:Tate_fully_faithful} and asserts that the canonical realization functor 
\begin{equation*}
    \Tate_{\aleph_0,R}^\mathrm{b} \longrightarrow \mathscr{D}(\cond_R)
\end{equation*}
from bounded countable Tate objects in perfect complexes  to the derived $\infty$-category $\zD(\cond_R)$ of condensed $R$-modules is fully faithful. 

Restricting to Tate objects that are \emph{bounded} (in the sense of having uniformly finite tor-amplitude, see Definition \ref{def:Tate_oo-cat}) is necessary, as we demonstrate in Section \ref{section:counterexample} by providing an explicit counterexample to full faithfulness for the $\infty$-category $\Tate_{\ale}(\Perf_R)$.\\

After a first version of this paper was released, P. Scholze kindly pointed out a mistake in the $\infty$-categorical full faithfulness: the property of $R$ being a finitely generated $\mathbb{Z}$-algebra is necessary (see Remark \ref{rem:scholze_counterexample}). He also made us aware of the recent paper \cite{Bergfalk_Lambie-Hanson}. There is some overlap between \cite[Section 2.2]{Bergfalk_Lambie-Hanson} and our results in Section \ref{infinity}, see Remark \ref{rem:proof_of_BL} for more details.
\\

\noindent{\textbf{Further directions}}. A natural continuation of this work, presently under investigation, would consist in defining 
\begin{itemize}
    \item a suitable first iteration $\mathsf{2-Solid}_R$ of the category $\mathsf{Solid}_R$ of solid $R$-modules;
    \item a fully faithful functor $$\mathsf{2-TC}_R\colon \mathsf{2-Tate}_{{\aleph}_0}(\mathsf{Proj}^{\mathrm{fg}}_R) \longhookrightarrow \mathsf{2-Solid}_R$$
    \item a solid version $$- \overrightarrow{\otimes}^{\blacksquare}- \colon \mathsf{Solid}_R \times \mathsf{Solid}_R \longrightarrow \mathsf{2-Solid}_R$$ of the normally oriented tensor product $$- \overrightarrow{\otimes} -\colon \mathsf{Tate}_{\aleph_0,R} \times \mathsf{Tate}_{\aleph_0,R} \longrightarrow \mathsf{2-Tate}_{\aleph_0}(\mathsf{Proj}^{\mathrm{fg}}_R)$$ 
        defined in \cite{BGHW},
\end{itemize}
\noindent in such a way that the diagram of 1-categories 
$$
\begin{tikzcd}[column sep = huge]
    \mathsf{Tate}_{\aleph_0,R} \times \mathsf{Tate}_{\aleph_0,R} \ar[d,"- \overrightarrow{\otimes} -"] \ar[r,"\mathsf{TC}_R \times \mathsf{TC}_R"] 
        & \mathsf{Solid}_R \times \mathsf{Solid}_R \ar[d,"- \overrightarrow{\otimes}^{\blacksquare} -"] \\
        \mathsf{2-Tate}_{\aleph_0}(\mathsf{Proj}^{\mathrm{fg}}_R) \ar[r,"{\mathsf{2-TC}_R}"] 
        & \mathsf{2-Solid}_R 
\end{tikzcd}
$$ 
commutes. The definition of $\mathsf{2-Solid}_R$ has been recently given by Vieri Sbandati in his Bachelor Thesis \cite{Vieri}, where he also proves that $\mathsf{2-Solid}_R$ is an Abelian full subcategory of $$\mathsf{2-Cond}_R:=\mathsf{Cond}(\mathsf{Cond}(\mathsf{Mod}_R))$$ stable under limits and colimits, and that this inclusion has a left adjoint (the so-called $2$-solidification functor). Once the 1-categorical situation is understood, a natural further step would be to investigate its $\infty$-categorical version.

Another interesting problem, that we were unable to solve, is to use our fully faithful functor $\mathsf{TC}_R$ in order to deduce Drinfeld's flat descent of Tate $R$-modules from the known flat descent of solid $R$-modules. K. Česnavičius kindly informed us of results that go in this direction. 

Finally, let us mention the recent approach of \emph{ultrasolid modules} proposed in \cite{Aparicio_ultrasolid}, which provides, when working over a field, a variant of solid modules using pro-modules. Although we do not pursue this here, it would be interesting to make a precise comparison between their theory and Tate modules.
\\

\noindent{\textbf{Acknowledgments}}. We started this project together with Andrea Maffei, and he contributed some important ideas; in particular, Section \ref{sec:1-cat_embedding} should be considered as joint work with him.
A first version of Proposition \ref{ff1} was obtained by Francesco Iacca in his Master Thesis \cite{Iacca}. We thank Chris Brav, Benjamin Hennion, Francesco Iacca for extremely useful discussions on the topics of this paper. We thank Peter Scholze for pointing out an error in the first version of this paper and for making us aware of the paper \cite{Bergfalk_Lambie-Hanson}.

\section{Embedding of Tate modules in condensed modules}
\label{sec:1-cat_embedding}

In this section we will produce, for any commutative (discrete) ring $R$, a fully faithful exact embedding of the exact category $\mathsf{Tate}_{{\aleph}_0,R}$  of countable Tate objects in the exact category $\mathsf{Proj}^{\mathrm{fg}}_{R}$ of finitely generated projective $R$-modules (as defined in \cite[Def. 5.23]{BGWTateExact}), into the abelian category $\mathsf{Cond}(\mathsf{Mod}_R)$ of condensed $R$-modules (as defined in \cite[Lectures 1 and 2]{CScond}).\\

We will write $\mathsf{Cond}(\mathcal{C})$ for the category of $\mathcal{C}$-valued sheaves on the site of profinite sets (endowed with the finitary jointly surjective topology), that are left Kan extended from their restriction to $\kappa$-small profinite sets for some uncountable strong limit cardinal $\kappa$ (see \cite[Def. 2.11]{CScond}). When $\mathcal{C}= \mathsf{Sets}$, we write $\mathsf{Cond}:= \mathsf{Cond}(\mathsf{Sets})$. When $R$ is a commutative ring and $\mathcal{C}= \mathsf{Mod}_R$, we will simply write $\mathsf{Cond}_R :=  \mathsf{Cond}(\mathsf{Mod}_{R})$.

We will write $\mathsf{TopMod}_R$ for the category of linearly topologized topological modules over the (discrete) ring $R$, and continuous $R$-linear morphisms.% {\color{blue}(or, if needed somewhere, its full subcategory of complete, separated, topological modules).}

\subsection{The embedding using topology}
We mostly follow the conventions and notations of \cite{BGWTateExact},  \cite{CScond}, and \cite{Drinfinite}.

\begin{defin} Let $R$ be a commutative ring. We let 
    \begin{itemize}
        \item $\mathsf{Tate}^{\mathrm{Dr}, \mathrm{el}}_{{\aleph}_0}(R)$ be the full subcategory of $\mathsf{TopMod}_{R}$ consisting of topological $R$-modules of the form $P\oplus Q^{\vee}$ where $P$ and $Q$ are discrete, countably generated projective $R$-modules, and $(-)^{\vee}$ denotes the topological linear dual.
        \item $\mathsf{Tate}^{\mathrm{Dr}}_{{\aleph}_0}(R)$ be the idempotent completion of $\mathsf{Tate}^{\mathrm{Dr}, \mathrm{el}}_{{\aleph}_0}(R)$ i.e. the full subcategory of $\mathsf{TopMod}_{R}$ of topological modules that are direct summands of objects in $\mathsf{Tate}^{\mathrm{Dr}, \mathrm{el}}_{{\aleph}_0}(R)$ (\cite[Def. 3.2.1]{Drinfinite}).
        \item $\mathsf{TopMod}^{\mathrm{cgwh}}_{R}$ be the full subcategory of $\mathsf{TopMod}_{R}$ consisting of topological modules whose underlying topological space is compactly generated and weakly Hausdorff.
    \end{itemize}
\end{defin}

If $\mathsf{Proj}^{\mathrm{fg}}_R$ is the category of finitely generated projective $R$-modules, we denote by $\mathsf{Tate}_{\aleph_0,R}:=\mathsf{Tate}_{\aleph_0}(\mathsf{Proj}^{\mathrm{fg}}_R)$ the 1-category of Tate $R$-modules whose size is countable, in the sense of \cite[Def. 5.3 and Def. 5.23]{BGWTateExact}. Recall from \cite{BGWTateExact} that there is a canonical functor 
\begin{equation}
    \label{eqn:realization_functor}
    \tau\colon \mathsf{Tate}_{{\aleph}_0,R} \longrightarrow \mathsf{TopMod}_{R}
\end{equation} 
obtained by realizing formal limits and colimits of Tate objects as actual limits and colimits in topological modules.

\begin{prop}\label{Drcg} Let $R$ be a commutative ring. The fully faithful inclusion $\mathsf{Tate}^{\mathrm{Dr}}_{{\aleph}_0}(R) \to \mathsf{TopMod}_{R}$ factors through the fully faithful inclusion $\mathsf{TopMod}^{\mathrm{cgwh}}_{R} \to \mathsf{TopMod}_{R}$, i.e. the underlying topological space of an object in $\mathsf{Tate}^{\mathrm{Dr}}_{{\aleph}_0}(R)$ is compactly generated and weakly Hausdorff.
\end{prop}

\begin{proof}
    By the proof of \cite[Thm. 5.26]{BGWTateExact}, the realization functor \eqref{eqn:realization_functor}  induces an exact equivalence of exact categories $$\mathsf{Tate}^{\mathrm{el}}_{{\aleph}_0}(\mathsf{Proj}^{\mathrm{fg}}_R) \simeq \mathsf{Tate}^{\mathrm{Dr}, \mathrm{el}}_{{\aleph}_0}(R).$$
    Therefore, \cite[Ex. 5.22]{BGWTateExact} tells us that  $\mathsf{Tate}^{\mathrm{Dr}, \mathrm{el}}_{{\aleph}_0}(R)$ is a split exact category, and any of its objects is a direct summand of $R(\!(t)\!)$. Since the underlying topological space of $R(\!(t)\!)$ is metrizable\footnote{Valuation theory gives the usual choice for a metric: $d(f,g):=2^{-\mathrm{ord}(f-g)}$.}, and any subspace of a metrizable space is metrizable, we get that any object in 
    $\mathsf{Tate}^{\mathrm{Dr}, \mathrm{el}}_{{\aleph}_0}(R)$ has an underlying metrizable topological space. The same argument shows that also any object in $\mathsf{Tate}^{\mathrm{Dr}}_{{\aleph}_0}(R)$ (which is, by definition, a direct summand of some object in  $\mathsf{Tate}^{\mathrm{Dr}, \mathrm{el}}_{{\aleph}_0}(R)$) is metrizable. We conclude since any metrizable topological space is compactly generated and (weakly) Hausdorff.
\end{proof}

\begin{cor}\label{fftotop} 
    The realization functor $\tau\colon \mathsf{Tate}_{{\aleph}_0,R} \to \mathsf{TopMod}_{R}$ is fully faithful and factors through the inclusion $\mathsf{TopMod}^{\mathrm{cgwh}}_{R} \hookrightarrow \mathsf{TopMod}_{R}$.
\end{cor}

\begin{proof} By \cite[Thm. 5.26]{BGWTateExact}, $\tau$ induces an exact equivalence between (split) exact categories $$\mathsf{Tate}_{{\aleph}_0,R} \simeq \mathsf{Tate}^{\mathrm{Dr}}_{{\aleph}_0}(R),$$ and we conclude by Proposition \ref{Drcg}.
\end{proof}

We consider the functor\footnote{This functor can be defined more generally for any T1 space, but one \emph{cannot} extend it to the whole category $\mathsf{TopMod}_R$, as explained in \cite[Warning 2.14]{CScond}.} 
\begin{equation*}\label{underlinemod} \underline{(-)} \colon\mathsf{TopMod}^{\mathrm{cgwh}}_{R} \longrightarrow \mathsf{Cond}_R \,\, \colon \, \, M \longmapsto \underline{M}:= C^{0}(\bullet,M).\end{equation*}
To be precise, the underlying condensed set of $\underline{M}$ is $C^0(-,U_{\mathrm{top}}(M))$, where $U_{\mathrm{top}}\colon \mathsf{TopMod}_{R} \to \mathsf{Top}$ is the obvious forgetful functor\footnote{This would be written $C^0(-,U_{\mathrm{top}}(M))= \underline{U_{\mathrm{top}}(M)}$ in the notations of \cite{CScond}, but we have not adopted this notation that would obviously be confusing here.}, while, for any $S\in \mathsf{ProFin}$, the $R$-module structure on $C^0(S,U_{\mathrm{top}}(M))$ is given by the obvious maps $$\xymatrix{C^0(S,U_{\mathrm{top}}(M)) \times C^0(S,U_{\mathrm{top}}(M)) \ar[r]^-{\sim} &  C^0(S,U_{\mathrm{top}}(M) \times U_{\mathrm{top}}(M)) \ar[rr]^-{C^0(S,\mathrm{sum}_M)} && C^0(S,U_{\mathrm{top}}(M))}$$

$$\xymatrix{C^0(S,R) \times C^0(S,U_{\mathrm{top}}(M)) \ar[r]^-{\sim} &  C^0(S,R \times U_{\mathrm{top}}(M)) \ar[rr]^-{C^0(S,\mathrm{mult}_M)} && C^0(S,U_{\mathrm{top}}(M))}$$ where $R$ is endowed with the discrete topology.

\begin{prop} The functor
    ${\underline{(-)}}\colon \mathsf{TopMod}^{\mathrm{cgwh}}_{R}  \to \mathsf{Cond}_R  $ is fully faithful.
\end{prop}

\begin{proof}  Consider the  functor $$C^0(-,\bullet)\colon \mathsf{Top}^{\mathrm{cgwh}} \longrightarrow \mathsf{Cond}\,\, \colon \,\, T \mapsto C^0(-,T),$$ and the commutative diagram
    \begin{equation}
        \label{eqn:diag_top_cond}
        \xymatrix{\mathsf{TopMod}^{\mathrm{cgwh}}_{R} \ar[r]^-{\underline{(-)}} \ar[d]_-{U_{\mathrm{Top}}} & \mathsf{Cond}_R \ar[d]^-{U_{\mathrm{Cond}}} \\
            \mathsf{Top}^{\mathrm{cgwh}} \ar[r]_{C^0(-, \bullet)} & \mathsf{Cond}
    }\end{equation}
    %Let $M, N \in \mathsf{TopMod}^{\mathrm{cgwh}}_{R}$.
    Recall from \cite[Prop. 1.7, Th. 2.16]{CScond} that the functor $C^0(-, \bullet)$ is fully faithful and factors through the full subcategory $\mathsf{qsCond}\subset \mathsf{Cond}$ of quasi-separated condensed sets, and that its corestriction admits a left adjoint $\mathsf{qsCond}\to \mathsf{Top}^{\mathrm{cgwh}}_R$ given by $X\mapsto X(*)_\mathrm{top}$. 
    In particular, given $M, N \in \mathsf{TopMod}^{\mathrm{cgwh}}_{R}$, the functor $C^0(-, \bullet)$ induces a bijection
    \begin{equation}
        \label{eqn:bij_Hom_Top_Cond}
        \mathrm{Hom}_{\mathsf{Top}^{\mathrm{cgwh}}}(U_{\mathrm{top}}(M),U_{\mathrm{top}}(N))
        \cong 
        \mathrm{Hom}_{\mathsf{Cond}}(U_{\mathrm{cond}}(\underline{M}),U_{\mathrm{cond}}(\underline{N}))
    \end{equation}
    with inverse $(-)(*)_{\mathrm{top}}$.
    Since $U_{\mathrm{Top}}$, $U_{\mathrm{Cond}}$ are faithful, commutativity of diagram \eqref{eqn:diag_top_cond} immediately implies that $\underline{(-)}$ is faithful. 
    To prove that it is also full, consider a morphism $f\colon\underline{M} \to \underline{N}$ in $\mathsf{Cond}_R$; we will show that it is of the form $\underline{g}$ for some map $g\colon M\to N$ in $\mathsf{TopMod}^{\mathrm{cgwh}}$.
    Using bijection \eqref{eqn:bij_Hom_Top_Cond} and commutative diagram \eqref{eqn:diag_top_cond} it suffices to prove that the continuous map $(U_{\mathrm{Cond}}f)(*)_{\mathrm{top}} \colon U_{\mathrm{top}}(M) \to U_{\mathrm{top}}(N)$ induced by $f$ is a morphism of topological $R$-modules, for then we could choose $g$ to be this morphism.
    Forgetting the topologies, we see that the map of sets $(U_{\mathrm{Cond}}f)(*)\colon U_{\mathrm{Set}}(M)\to U_{\mathrm{Set}}(N)$ is the one underlying the $R$-modules morphism $f(*)\colon \underline{M}(*)\to \underline{N}(*)$, which gives the desired result.
\end{proof}

The previous two results immediately imply the following proposition.
\begin{prop}\label{ff1} The composite functor
    $$\mathrm{TC_R}\colon\xymatrix{\mathsf{Tate}_{{\aleph}_0 , R} \ar[r]^-{\tau} & \mathsf{TopMod}^{\mathrm{cgwh}}_{R} \ar[r]^-{\underline{(-)}} & \mathsf{Cond}_R  }$$ is fully faithful.
\end{prop}

\subsection{A non-topological description of the embedding}
We may give a reformulation of the functor $\mathrm{TC}_R$ of Prop. \ref{ff1} that avoids going through topological modules.\\

We start by denoting 
\begin{equation}
    \label{eqn:def_functor_a(const)}
    a(\mathrm{const}_{-})\colon\mathsf{Mod}_R \longrightarrow \mathsf{Cond}_R
\end{equation} 
the functor sending a module $P$ to the sheaf (condensed $R$-module) associated to the constant presheaf $\mathrm{const}_P$ with value $P$ on any $S \in \mathsf{ProFin}$. The functor $a(\mathrm{const}_{-})$ is fully faithful and left adjoint to the functor given by evaluation of a condensed $R$-module at the singleton $* \in \mathsf{ProFin}$. \\
Let $(-)_{\delta}\colon\mathsf{Mod}_R \to \mathsf{TopMod}_R$ the functor sending $P$ to itself endowed with the discrete topology $P_{\delta}$ (note that $R$ is always tacitly endowed with the discrete topology, so $P_{\delta}$ is indeed a topological $R$-module). Since discrete topological spaces are compactly generated and (weakly) Hausdorff, we actually have a functor $(-)_{\delta}\colon\mathsf{Mod}_R \to \mathsf{TopMod}^{\mathrm{cgwh}}_R$ which is a left adjoint.

\begin{lem}\label{Andrea} There is a natural isomorphism 
    $$\alpha\colon a(\mathrm{const}_{-}) \longrightarrow \underline{(-)} \circ (-)_{\delta}$$
    of functors $\mathsf{Mod}_R \longrightarrow \mathsf{Cond}_R$.
\end{lem}
\begin{proof}
    Since $a(\mathrm{const}_{-})$ is the fully faithful left adjoint to evaluation at $*$, we define $\alpha$ as the adjunct of the functorial isomorphism $ P \simeq \underline{(P_\delta )}(*)= C^0(*, P_\delta)$ in $\mathsf{Mod}_R$. 
    Let us prove that, for any $P\in \mathsf{Mod}_R$, $\alpha(P)\colon a(\mathrm{const}_{P}) \to \underline{P_\delta  } $ is indeed an isomorphism in $\mathsf{Cond}_R$. Since $\alpha(P)$ is a morphism in $\mathsf{Cond}_R$, and the forgetful functor $ U_{\mathrm{Cond}}\colon\mathsf{Cond}_R \to \mathsf{Cond}$ is conservative, it is enough to prove that $U_{\mathrm{Cond}}(\alpha (P))$ is an isomorphism in $\mathsf{Cond}$.  Let $K$ be an arbitrary set, and consider $\alpha^{\mathrm{pre}}(K)\colon\mathrm{const}_K \to C^0 (-, K_{\delta})$ the morphism in $\mathsf{Fun}(\mathsf{ProFin}^{op}, \mathsf{Set})$, defined by the inclusion $\alpha^{\mathrm{pre}}(K)(S) \colon K \to C^0(S, K_\delta)$ of constant functions into locally constant functions, for any $S \in \mathsf{ProFin}$. It will be enough to show that the morphism $\alpha^{\mathrm{pre}}(K)$ exhibits $C^0 (-, K_{\delta})$ as the sheafification of the constant presheaf $\mathrm{const}_K$. In order to achieve this, let $X \in \mathsf{Cond}$, and $\varphi\colon\mathrm{const}_K \to X$ a morphism in $\mathsf{Fun}(\mathsf{ProFin}^{op}, \mathsf{Set})$: we will produce a unique $\psi\colon C^0(-, K_\delta) \to X$ such that $\psi \circ \alpha^{\mathrm{pre}}(K)= \varphi$. Let $S \in \mathsf{ProFin}$, and $f \in C^0(S,K_\delta)$. For $k \in K$, let $S_k := f^{-1}(k)$ which is open in $S$, and $\cup_{k\in K} S_k =S$. Since $S$ is quasicompact, there is a finite subset $K' \subset K$ such that $\cup_{k'\in K'} S_{k'} =S$
    (a finite open subcover), and the union is disjoint. But $X$ is a sheaf, hence $X(S)=X(\coprod_{k' \in K'} S_{k'}) \simeq \prod_{k' \in K'}X(S_{k'})$, and we define $\psi(S)\colon C^0(S, K_\delta) \to X(S)$ by sending $f$ to the family $(\varphi(S_{k'})(k))_{k' \in K'}$. It is then easy to show that such a $\psi$ is the unique map $\psi\colon C^0(-, K_\delta) \to X$ such that $\psi \circ \alpha^{\mathrm{pre}}(K)= \varphi$.

\end{proof}

Now we want to extend the functor $a(\mathrm{const}_{-})$ to a functor $$\mathrm{TC}'_R\colon \mathsf{Tate}_{{\aleph}_0 , R} \longrightarrow \mathsf{Cond}_R$$ defined on objects in $\mathsf{Tate}^{\textrm{el}}_{{\aleph}_0 , R} \subset \mathsf{Tate}_{{\aleph}_0 , R}$ by
\begin{equation}\label{indproconst} \mathsf{Tate}^{\textrm{el}}_{{\aleph}_0 , R} \ni \text{``$\mathrm{lim}_i \,\mathrm{colim}_j$''}  \, P_{ij} \longmapsto \mathrm{lim}_i \,\mathrm{colim}_j \, a(\mathrm{const}_{P_{ij}}) \, \in \mathsf{Cond}_R
\end{equation}
where the colimit and limit in the target is taken in $\mathsf{Cond}_R$.\\
The proof of the next Proposition shows how to construct $\mathrm{TC}'_R$ rigorously.

\begin{prop}\label{TC'} The functor $a(\mathrm{const}_{-})\colon \mathsf{Proj}^{\mathrm{fg}}_R \longrightarrow \mathsf{Cond}_R$ extends to an exact functor $$\mathrm{TC}'_R\colon \mathsf{Tate}_{{\aleph}_0 , R} \longrightarrow \mathsf{Cond}_R$$ defined on objects as in (\ref{indproconst}).
\end{prop}

\begin{proof} Since $\mathsf{Cond}_R$ is (Abelian hence) idempotent complete, it will be enough to construct an extension $\mathsf{Tate}^{\textrm{el}}_{{\aleph}_0 , R} \to \mathsf{Cond}_R$. By  \cite[Thm. 5.4 and Prop. 5.20]{BGWTateExact} $\mathsf{Tate}^{\textrm{el}}_{{\aleph}_0 , R}$ is the smallest subcategory of $\mathsf{Ind}^{\mathrm{a}}_{{\aleph}_0}\mathsf{Pro}^{\mathrm{a}}_{{\aleph}_0}(\mathsf{Proj}^{\mathrm{fg}}_R)$ containing both $\mathsf{Pro}^{\mathrm{a}}_{{\aleph}_0}(\mathsf{Proj}^{\mathrm{fg}}_R)$ and $\mathsf{Ind}^{\mathrm{a}}_{{\aleph}_0}(\mathsf{Proj}^{\mathrm{fg}}_R)$, and closed under extensions, and it is split exact. Therefore it suffices to produce exact functors
    $$F^{\mathrm{Pro}}\colon \mathsf{Pro}^{\mathrm{a}}_{{\aleph}_0}(\mathsf{Proj}^{\mathrm{fg}}_R) \to \mathsf{Cond}_R \qquad F^{\mathrm{Ind}}\colon \mathsf{Ind}^{\mathrm{a}}_{{\aleph}_0}(\mathsf{Proj}^{\mathrm{fg}}_R) \to \mathsf{Cond}_R$$ both restricting to  $$a(\mathrm{const}_{-}) \colon \mathsf{Proj}^{\mathrm{fg}}_R \longrightarrow \mathsf{Cond}_R $$ 
    Now, $\mathsf{Cond}_R$ is complete and co-complete (i.e. satisfies (AB3) and (AB3*)\cite[Thm. 2.2]{CScond}), in particular it has $\aleph_0$-filtered colimits and $\aleph_0$-cofiltered limits, so by the universal property of $\mathsf{Pro}^{\mathrm{a}}_{{\aleph}_0}$
and $\mathsf{Ind}^{\mathrm{a}}_{{\aleph}_0}$ we get $F^{\mathrm{Pro}}$, and $F^{\mathrm{Ind}}$ with the desired properties. Since the exact structure on $\mathsf{Tate}_{\aleph_0}(\mathsf{Proj}^{\mathrm{fg}})$ is $\mathsf{Tate}_{\aleph_0}(\mathsf{Ex^{\mathsf{Proj}^{\mathrm{fg}}_R}})$, and the exact structure $\mathsf{Ex}^{\mathsf{Proj}^{\mathrm{fg}}_R}$ in $\mathsf{Proj}^{\mathrm{fg}}_R$ is split, $\mathrm{TC}'_R$ is, by definition, exact.   \end{proof}

\begin{prop}\label{ff2} There is an isomorphism $\theta\colon  \mathrm{TC}_R \simeq \mathrm{TC}'_R$ of functors $\mathsf{Tate}_{{\aleph}_0 , R} \longrightarrow \mathsf{Cond}_R$. In particular, $\mathrm{TC}'_R$ is fully faithful since $\mathrm{TC}_R$ is (Prop. \ref{ff1}), and both are exact, since $\mathrm{TC}'_R$ is (Prop. \ref{TC'}).
\end{prop}

\begin{proof} It is enough to define a natural isomorphism $$\theta\, |_{ \mathsf{Tate}^{\mathrm{el}}_{{\aleph}_0 , R}}\colon  \mathrm{TC}'_R \, |_{ \mathsf{Tate}^{\mathrm{el}}_{{\aleph}_0 , R}}  \longrightarrow  \mathrm{TC}_R \, |_{\mathsf{Tate}^{\mathrm{el}}_{{\aleph}_0 , R}} $$ as follows. Using that $a(\mathrm{const}_{-})$ is left adjoint together with Lemma \ref{Andrea}, if  $\text{``$\mathrm{lim}_i \,\mathrm{colim}_j$''}  \, P_{ij} \in \mathsf{Tate}^{\mathrm{el}}_{{\aleph}_0 , R}$, we have 
    \begin{align*}
        \mathrm{TC}'_R (\text{``$\mathrm{lim}_i \,\mathrm{colim}_j$''}  \, P_{ij}) 
        & = \mathrm{lim}_i \, \mathrm{colim}_j \, a(\mathrm{const}_{P_{ij}}) \\
        & \simeq \mathrm{lim}_i \, a(\mathrm{const}_{\mathrm{colim}_j \, P_{ij}}) \\
        & \simeq \mathrm{lim}_i \, \underline{(\mathrm{colim}_j \, P_{ij})_\delta} \\
        & = \mathrm{lim}_i \, C^0(-,(\mathrm{colim}_j \, P_{ij})_\delta ).
    \end{align*}
    Since the discrete topology functor $(-)_\delta \colon \mathsf{Mod}_R \to \mathsf{TopMod}^{\mathrm{cgwh}}_{R}$ is left adjoint (to the forgetful one), %{\color{purple} and the forgetful functor $$U_{\mathrm{top}}\colon \mathsf{TopMod}^{\mathrm{cgwh}}_{R} \to \mathsf{Top}$$ is right adjoint and commutes with $\mathbb{N}$-filtered colimits along (closed) subspace immersions [I don't see where this last fact plays any role, I think the fact that $(-)_\delta$ is a left adjoint is enough], }
    we finally get
    $$\mathrm{lim}_i \, C^0(-,(\mathrm{colim}_j \, P_{ij})_\delta ) \simeq \mathrm{lim}_i \, C^0(-,\mathrm{colim}_j \, (P_{ij})_\delta )\simeq C^0(-,\mathrm{lim}_i \, \mathrm{colim}_j \, (P_{ij})_\delta ) = \mathrm{TC}_R (\text{``$\mathrm{lim}_i \,\mathrm{colim}_j$''}  \, P_{ij}). $$
    %{\color{blue} [Doubt: is the second $\simeq$ correct??]} {\color{purple}[Hugo: I think it is fine, since $\underline{(-)}$ is a right adjoint.]}
\end{proof}

Recall from \cite[Proposition 3.16]{Andreychev_2024} that for any ring $R$ can be endowed with an analytic ring structure $R_\blacksquare$ by considering the functor of measures $R_\blacksquare[-]$ that sends a profinite set $S\cong\lim_i S_i$ to the condensed $R$-module
\[
    R_\blacksquare[S] := \colim_{R' \subseteq R} \lim_i R'[S_i],
\]
where the colimit runs over finitely generated subrings $R'$ of $R$. The associated category of solid modules, defined in \cite[Proposition 7.5]{CScond}, will be denoted  $\mathsf{Solid}_R$.\footnote{In \cite{CScond}, this category of solid $R$-modules is denoted $\mathsf{Mod}_{R_\blacksquare}^\mathrm{cond}$.}

Since the condensed $R$-module $a(\mathrm{const}_{P})$ associated to any discrete $R$-module $P$ is solid, and since the full subcategory  $\mathsf{Solid}_R \subset \mathsf{Cond}_R$ is closed under colimits and limits, we get that $\mathrm{TC}_R$ and $\mathrm{TC}'_R$ factor via fully faithful embeddings (denoted with the same symbols) 
\begin{equation*}
    \label{eqn:functor_TS}
    \mathrm{TS}_R \simeq \mathrm{TS}'_R \colon\mathsf{Tate}_{{\aleph}_0 , R} \longrightarrow \mathsf{Solid}_R.
\end{equation*}

\begin{rem}
    Alternatively, one could consider the full subcategory $\mathsf{Mod}_R(\mathsf{Solid}_\mathbb{Z})\subseteq \mathsf{Cond}_R$ consisting of those condensed $R$-modules whose underlying condensed abelian group is solid; in other words, this amounts to replacing $R_\blacksquare$ by the analytic ring $(R,\mathbb{Z})_\blacksquare$. The same argument as before shows that the embedding $\mathrm{TC}_R$ also factors through this full subcategory.
    % {\color{blue}[Here we should be clear about which of the 2 possible definitions of $\mathsf{Solid}_R$ we use (in order the Remark to be true). Which one? Both? Need $R$ to be finitely generated over $\mathbb{Z}$?]}
    % {\color{purple}[Hugo: It think it works for both definitions: in both cases $R$ is solid and $\mathsf{Solid}_R$ is stable under colimits so that all of $\mathsf{Mod}_R$ belongs to $\mathsf{Solid}_R$.]}
    % \label{goestosolid} 
    % For any (discrete) $R$-module $P$, the condensed $R$-module $a(\mathrm{const}_{P})$ is a solid $R$-module. Since the abelian full subcategory  $\mathsf{Solid}_R \subset \mathsf{Cond}_R$ of solid $R$-modules is closed under colimits and limits, we get that $\mathrm{TC}_R$ and $\mathrm{TC}'_R$ factor via fully faithful embeddings (denoted with the same symbols) 
    % \begin{equation*}
    %     \label{eqn:functor_TS}
    %     \mathrm{TS}_R \simeq \mathrm{TS}'_R \colon\mathsf{Tate}_{{\aleph}_0 , R} \longrightarrow \mathsf{Solid}_R.
    % \end{equation*}
\end{rem}

\section{Essential image of the embedding over a ring of finite type}
\label{sec:essential_image}

We start with a first description of the essential image of the embedding $\mathrm{TS}$, which holds for any base ring.
\begin{prop}
    The essential image of $\mathrm{TS}_R$ consists of direct summands of $R(\!(t)\!)= R^\mathbb{N}\oplus R^{(\mathbb{N})}$.
\end{prop}

\begin{proof} Observe that any object in $\mathsf{Tate}_{{\aleph}_0,\, R}$ is a split idempotent of $R(\!(t)\!)$. Then, $\mathrm{TS}_R\colon\mathsf{Tate}_{{\aleph}_0,\, R}  \longrightarrow \mathsf{Solid}_R $ being fully faithful between idempotent complete categories, it induces an equivalence between the full subcategories of (split) idempotents of $R(\!(t)\!)$ and of $\mathrm{TS}_R(R(\!(t)\!))$.
\end{proof}

Until the end of this section, we will assume that the base ring $R$ is a finitely generated $\mathbb{Z}$-algebra. 
In this case, we will show the following characterization of the essential image of Tate modules inside solid ones.

\begin{thm}
    \label{thm:ess_image}
    For $R$ a finitely generated $\mathbb{Z}$-algebra, the essential image of the functor 
    $\mathrm{TS}_R \colon\mathsf{Tate}_{{\aleph}_0 , R} \to \mathsf{Solid}_R$
    is the smallest full subcategory of $\mathsf{Solid}_R$ that contains the free module $\bigoplus_{\mathbb{N}} R$ of infinite countable rank and is stable under direct sums, retracts and duals.
\end{thm}

The proof will be given at the end of the section. Before that, we need to analyze linear duality in the categories $\mathsf{Tate}_{\aleph_0,R}$ and $\mathsf{Solid}_{R}$, and to compare them.

\subsection{Duality for Tate modules} By \cite[Proposition 3.5]{BGHW}, any exact equivalence $\Psi\colon \zC^{\mathrm{op}}\simeq \zC$ of an idempotent complete exact category $\zC$ extends to an exact equivalence 
\begin{equation}
    \label{eqn:duality_functor_Tate}
    \mathsf{Tate}^{\mathrm{el}}(\zC)^{\mathrm{op}} \stackrel{\sim}{\longrightarrow}
    \mathsf{Tate}^{\mathrm{el}}(\zC)
\end{equation}
that moreover exchanges Pro and Ind objects. 

We briefly recall the construction of this duality functor. First consider the functor 
\begin{equation*}
    \Phi \colon \mathsf{Tate}^{\mathrm{el}}(\zC) \longrightarrow \mathsf{Pro}^{\mathrm{a}}\mathsf{Ind}^{\mathrm{a}}(\zC)
\end{equation*}
sending an elementary Tate object $V$ to the Pro-object $L\mapsto V/L$ indexed by the poset of lattices\footnote{A \emph{lattice}  of $V$ is an admissible monic $L\to V$ such that $L\in \mathsf{Pro}^{\mathrm{a}}(\zC)$ and $V/L\in \mathsf{Ind}^{\mathrm{a}}(\zC)$, cf \cite[Definition 6.1.]{BGWTateExact}.} $L$ of $V$. 
The composite 
\begin{equation*}
    \mathsf{Tate}^{\mathrm{el}}(\zC)^{\mathrm{op}} \longrightarrow \left(\mathsf{Pro}^{\mathrm{a}}\mathsf{Ind}^{\mathrm{a}}(\zC)\right)^{\mathrm{op}} \cong \mathsf{Ind}^{\mathrm{a}}\mathsf{Pro}^{\mathrm{a}}\left(\zC^{\mathrm{op}}\right)
\end{equation*}
of $\Phi^{\mathrm{op}}$ with the canonical isomorphism factors through elementary Tate objects and restricts to an equivalence 
\begin{equation*}
    \overline{\Phi} \colon
    \mathsf{Tate}^{\mathrm{el}}(\zC)^{\mathrm{op}} \stackrel{\sim}{ \longrightarrow} 
    \mathsf{Tate}^{\mathrm{el}}(\zC^{\mathrm{op}}). 
\end{equation*}
The duality functor \eqref{eqn:duality_functor_Tate} is then obtained as the composite of $\overline{\Phi}$ with the functor  
$\mathsf{Tate}^{\mathrm{el}}(\Psi) \colon
\mathsf{Tate}^{\mathrm{el}}(\zC^{\mathrm{op}})
\stackrel{\sim}{ \longrightarrow} 
\mathsf{Tate}^{\mathrm{el}}(\zC)$
induced by $\Psi$.

Applying this construction to the linear duality functor $(\mathsf{Proj}^{\mathrm{fg}}_R)^\mathrm{op} \to \mathsf{Proj}^{\mathrm{fg}}_R $, we obtain an exact equivalence 
\begin{equation}
    \label{eqn:linear_duality_Tate}
    (-)^\vee \colon \mathsf{Tate}_{{\aleph}_0,R}^{\mathrm{op}} \stackrel{\sim}{ \longrightarrow} \mathsf{Tate}_{{\aleph}_0,R}.
\end{equation}
By the explicit description given above, one easily sees that there are canonical isomorphisms
\begin{equation}
    \label{eqn:computation_dual_Tate}
    \left( R^\mathbb{N} \right)^\vee \cong R^{(\mathbb N)}
    \qquad
    \text{and}
    \qquad
    \left( R^{(\mathbb N)} \right)^\vee \cong R^\mathbb{N}.
\end{equation}

\subsection{Duality for solid modules}
Recall that in this section we assume that the base ring $R$ is a finitely generated $\mathbb{Z}$-algebra. In particular, for any set $I$ the product $\prod_I R$ is a compact projective object of $\mathsf{Solid}_R$.
\begin{lem}
    \label{lem:duality_solid}
    The solid modules $\prod_\mathbb{N} R$ and $\bigoplus_\mathbb{N} R$ are linear dual to one another.
\end{lem}
\begin{proof}
    It is clear that $\left( \bigoplus_\mathbb{N} R \right)^\vee \cong \prod_\mathbb{N} R$. We now prove that $\left( \prod_\mathbb{N} R \right)^\vee \cong \bigoplus_\mathbb{N} R$.
    We will use that the solid modules $\prod_I R$ for varying sets $I$ form a family of compact projective generators of $\mathsf{Solid}_R$.
    For any such set $I$ we have bijections 
    \begin{align*}
        \Hom_{\mathsf{Solid_R}}\left(\prod_I R, \left( \prod_\mathbb{N} R \right) ^\vee\right) 
        & \cong \Hom_{\mathsf{Solid_R}}\left(\prod_I R \otimes^{\blacksquare}_R \prod_\mathbb{N} R, R\right) \\
        & \cong \Hom_{\mathsf{Solid_R}}\left(\prod_{I\times \mathbb{N}} R, R\right) \\
        & \cong \bigoplus_\mathbb{N}\Hom_{\mathsf{Solid_R}}\left(\prod_{I} R, R\right) \\
        & \cong \Hom_{\mathsf{Solid_R}}\left(\prod_{I} R, \bigoplus_\mathbb{N}R\right)
    \end{align*}
    where the second isomorphism follows from \cite[Proposition 6.3.]{CScond}, the third from fully faithfulness of $\mathrm{TS}$ (Proposition \ref{ff1}) and the last one from compactness of $\prod_I R$.
\end{proof}

\subsection{Characterization of the essential image}
%\label{sec:essential_image}

\begin{lem}
    \label{lem:duality_endomorphisms}
    Let $f$ be an endomorphism of $R^\mathbb{N} \oplus R^{(\mathbb{N})}$ in $ \mathsf{Tate}^{\mathrm{el}}_{\aleph_0,R}$. Then $\mathrm{TS}(f^\vee) = \mathrm{TS}(f)^\vee$.
\end{lem}
\begin{proof}
    We can decompose $f$ as a matrix 
    $\begin{pmatrix} 
        f_{11} & f_{12} \\ 
        f_{21} & f_{22} 
    \end{pmatrix}$ of morphisms $f_{ij} \colon A_i \to A_j$, where we write $A_1 := \prod_\mathbb{N} R$ and $A_2 := \bigoplus_\mathbb{N} R$. It then suffices to show that $\mathrm{TS}(f_{ij}^\vee) = \mathrm{TS}(f_{ij})^\vee$ for every $i,j \in \{1,2\}$. 

    We start with the case of 
    $f_{11} \in \mathrm{End}_{\mathsf{Tate}_R%^{\mathrm{el}}_{\aleph_0,R}
    }(\prod_\mathbb{N} R)$.
    Observe that by definition of morphisms of Tate objects, there are isomorphisms of abelian groups
    \[
        {\End_{\mathsf{Tate}_R}\left(\prod_\mathbb{N}R\right)} \cong {\prod_\mathbb{N}\bigoplus_\mathbb{N}R} \cong \End_{\mathsf{Tate}_R}\left(\bigoplus_\mathbb{N}R\right).
    \]
    Since $\mathrm{TS}$ is fully faithful, we have similar isomorphisms for endomorphisms in solid $R$-modules ${\End_{\mathsf{Cond}_R}\left(\prod_\mathbb{N}R\right)} \cong {\prod_\mathbb{N}\bigoplus_\mathbb{N}R}
    \cong \End_{\mathsf{Cond}_R}\left(\bigoplus_\mathbb{N}R\right) .$ Using the isomorphisms \eqref{eqn:computation_dual_Tate} and the ones provided by Lemma \ref{lem:duality_solid}, we need to prove that in the diagram
    % https://q.uiver.app/#q=WzAsNSxbMCwwLCJcXEVuZF97XFxtYXRoc2Z7VGF0ZX1fUn0oXFxwcm9kX1xcbWF0aGJie059UikiXSxbMSwxLCJcXHByb2RfXFxtYXRoYmJ7Tn1cXGJpZ29wbHVzX1xcbWF0aGJie059UiJdLFsyLDAsIlxcRW5kX3tcXG1hdGhzZntUYXRlfV9SfShcXGJpZ29wbHVzX1xcbWF0aGJie059UikiXSxbMCwyLCJcXEVuZF97XFxtYXRoc2Z7Q29uZH1fUn0oXFxwcm9kX1xcbWF0aGJie059UikiXSxbMiwyLCJcXEVuZF97XFxtYXRoc2Z7Q29uZH1fUn0oXFxiaWdvcGx1c19cXG1hdGhiYntOfVIpIl0sWzAsMV0sWzAsMiwiKC0pXlxcdmVlIl0sWzIsMV0sWzAsM10sWzMsMV0sWzMsNCwiKC0pXlxcdmVlIiwyXSxbNCwxXSxbMiw0XV0=
    \[\begin{tikzcd}
        {\End_{\mathsf{Tate}_R}(\prod_\mathbb{N}R)} && {\End_{\mathsf{Tate}_R}(\bigoplus_\mathbb{N}R)} \\
                                                    & {\prod_\mathbb{N}\bigoplus_\mathbb{N}R} \\
        {\End_{\mathsf{Cond}_R}(\prod_\mathbb{N}R)} && {\End_{\mathsf{Cond}_R}(\bigoplus_\mathbb{N}R)}
        \arrow["{(-)^\vee}", from=1-1, to=1-3]
        \arrow[from=1-1, to=2-2,"\cong" above]
        \arrow["\mathrm{TS}","\cong" left, from=1-1, to=3-1]
        \arrow[from=1-3, to=2-2,"\cong" above]
        \arrow["\mathrm{TS}","\cong" left, from=1-3, to=3-3]
        \arrow[from=3-1, to=2-2,"\cong" above]
        \arrow["{(-)^\vee}"', from=3-1, to=3-3]
        \arrow[from=3-3, to=2-2,"\cong" above]
    \end{tikzcd}\]
    the outer square commutes. To do so, it suffices to show that every inner triangle commutes. We already know that the left and right triangle commutes, by construction. Now for the top (respectively the bottom) triangle, commutativily follows from inspection of the explicit construction of the duality functor \eqref{eqn:linear_duality_Tate} (resp. that of the proof of Lemma \ref{lem:duality_solid}). This proves that $\mathrm{TS}(f_{11}^\vee) = \mathrm{TS}(f_{11})^\vee$. The case of an endomorphism $f_{22}$ of $\bigoplus_\mathbb{N} R$ is completely similar. 

    We now turn to the off-diagonal morphisms, but spell out the argument only for $f_{12}\colon \prod_\mathbb{N} R\to \bigoplus_\mathbb{N}R$, the case of $f_{21}$ being similar. In that case the corresponding abelian group of morphisms can be computed as  
    \begin{equation}
        \label{eqn:computation_f_12}
        {\Hom_{\mathsf{Cond}_R}\left(\prod_\mathbb{N}R, \bigoplus_\mathbb{N}R \right)} \stackrel{\mathrm{TC}}{\cong} 
        {\Hom_{\mathsf{Tate}_R}\left(\prod_\mathbb{N}R, \bigoplus_\mathbb{N}R \right)} \cong {\bigoplus_\mathbb{N}\bigoplus_\mathbb{N}R}.
    \end{equation}
    The dual morphism $f_{12}^\vee \colon \left( \bigoplus_\mathbb{N} R \right) ^\vee \to \left( \prod_\mathbb{N} R \right) ^\vee $ corresponds via the isomorphisms \eqref{eqn:computation_dual_Tate} to a map $\prod_\mathbb{N} R\to \bigoplus_\mathbb{N} R$, which can be seen to be $f_{12}$, using the explicit description of the linear duality functor in Tate modules. Similarly on the condensed side, tracing back through the isomorphisms \eqref{eqn:computation_f_12} and those in the proof of Lemma \ref{lem:duality_solid}, we see that $\mathrm{TS}(f_{12})^\vee$ can be identified with $\mathrm{TS}(f_{12})$, which gives the desired result. 
\end{proof}

\begin{prop}
    \label{preserveduals_bis} 
    The functor $\mathrm{TS}\colon\mathsf{Tate}_{{\aleph}_0,R} \longrightarrow \mathsf{Solid}_R  $ preserves duals.
\end{prop}
\begin{proof}
    Let $T$ be a countable Tate $R$-module. By \cite[Example 7.6.]{BGWTateExact} $T$ is a direct summand of $R(\!(t)\!)$. Let $p$ be an idempotent of $R(\!(t)\!)$ with kernel $T$. Then we have isomorphisms 
    \begin{equation*}
        \mathrm{TS}(T)^\vee \cong (\ker\mathrm{TS}(p))^\vee
        \cong \coker \left( \mathrm{TS}(p)^\vee \right) 
        \cong \coker \left( \mathrm{TS}(p^\vee) \right) 
        \cong \mathrm{TS}(T^\vee),
    \end{equation*}
    where the third one is given by Lemma \ref{lem:duality_endomorphisms}.
\end{proof}

We can now provide the characterization of the essential image of the embedding $ \mathsf{Tate}_{{\aleph}_0,R} \hookrightarrow\mathsf{Solid}_R$. 

\begin{proof}[Proof of Theorem \ref{thm:ess_image}]
    Let $\cC$ denote the smallest full subcategory of $\mathsf{Solid}_R$ that contains $\bigoplus_{\mathbb{N}} R$ and is stable under direct sums, retracts and duals.

    We first prove that the essential image $\im(\mathrm{TS})$ is contained in $\cC$. %Since, it suffices to show that $\mathrm{TS}(X) \in \cC$ for every countable elementary Tate object $X \in \mathsf{Tate}_{{\aleph}_0 , R}^{\mathrm{el}}$. 
    By \cite[Example 7.6]{BGWTateExact} every countable Tate object  is a summand of $R^\mathbb{N}\oplus R^{(\mathbb{N})}$; since  $\cC$ is idempotent complete, it suffices to show that the latter object is in $\cC$. By definition $\cC$ contains $\bigoplus_\mathbb{N} R$, hence its dual $\prod_\mathbb{N} R$, therefore also their sum $\bigoplus_\mathbb{N} R\oplus \prod_\mathbb{N} R$. 

    We now show the converse inclusion. First, it is clear that $\im(\mathrm{TS})$ contains $\bigoplus_\mathbb{N} R$ and is stable under direct sums. Since any retract is a kernel of a projector and  $\mathrm{TS}$ is fully faithful (by Proposition \ref{ff1}) and exact (by Proposition \ref{ff2}), it follows that $\im(\mathrm{TS})$ is also stable under retracts.
    Finally, stability under taking duals is given by Proposition \ref{preserveduals_bis}. This shows that $\cC \subseteq \im(\mathrm{TS})$, concluding the proof.
\end{proof}

\section{\texorpdfstring{$\infty$}{Infinity}-categorical embedding in the bounded case}
\label{infinity}

Let $R$ be a ring of finite type over $\mathbb{Z}$.
In this section, we show that countable Tate objects in perfect complexes of  $R$-modules, if uniformly bounded in the sense of tor-amplitude, embed fully faithfully into the derived $\infty$-category of condensed $R$-modules.

\begin{defin}
    \label{def:Tate_oo-cat}
    We will consider the following $\infty$-categories.
    \begin{itemize}
        \item 
            Following \cite[Definition 2.1]{hennion:tate}, we let $\Tate_{R} := \Tate(\Perf_R)$ be the $\infty$-category of Tate objects in perfect complexes of $R$-modules. Recall that this is the smallest stable idempotent complete subcategory of $\Pro\,\Ind(\Perf_R)$ containing both the essential images of $\Pro(\Perf_R)$ and $\Ind(\Perf_R)\simeq \zD(R)$.
        \item Restricting to countable cofiltered  diagrams gives the full subcategory $\Pro_\ale(-)\subset \Pro(-)$  of \emph{countable} pro-objects; and dually for $\Ind_\ale(-)\subset \Ind(-)$.
            We define the  $\infty$-category $\Tate_{\aleph_0,R}$ of \emph{countable Tate objects} in $\Perf_R$ as the intersection of $\Tate_R$ with $\Pro_\ale(\Ind_\ale(\Perf_R))$.
            %, by which we mean those Tate objects that can be expressed using only \emph{countable} filtered and cofiltered diagrams.
        \item 
            Given integers $a$ and $b$, let $\Perf_R^{[a,b]}$ denote the full subcategory of the stable $\infty$-category $\Perf_R$ of perfect complexes of $R$-modules consisting of objects with tor-amplitude in $[a,b]$, where we use cohomological grading conventions. 
             We may form the $\infty$-category $\Pro_\ale \Ind_\ale \left( \Perf_R^{[a,b]} \right) $ of countable pro-ind-perfect complexes that are \emph{uniformly of tor-amplitude in $[a,b]$}. The intersection of this $\infty$-category with $\Tate_{\aleph_0,R}$ will be denoted $\Tate_{\aleph_0,R}^{[a,b]}$.
        \item Taking the union of these subcategories $\Tate_{\aleph_0,R}^{[a,b]}$ for all $a,b \in \mathbb{Z}$ gives the $\infty$-category of \emph{bounded countable Tate objects}, denoted $\Tate_{\ale,R}^\mathrm{b}$.
        \item Given an $\infty$-category $\zC$ with finite limits, we let $\cond(\zC)$ be the $\infty$-category of condensed objects in $\zC$, defined as accessible $\zC$-valued sheaves on profinite sets. For $\zC = \zD(R)$, note that there is an equivalence $\mathsf{Cond}(\zD(R)) \simeq \zD(\mathsf{Cond}_R)$ of stable $\infty$-categories.
    \end{itemize}
\end{defin}

Observe that the fully faithful functor $\mathsf{Mod}_R\hookrightarrow \mathsf{Cond}_R$ defined in \eqref{eqn:def_functor_a(const)} induces a fully faithful functor of stable $\infty$-categories $\zD(R)\hookrightarrow\zD(\mathsf{Cond}_R)$, which is  left adjoint to evaluation at $*\in \mathsf{ProFin}$, and factors through the stable full subcategory $\zD(\solid_R)$ of solid objects. Beware that, as opposed to the previous sections, we will not distinguish between an object in $\zD(R)$ and its image through this functor.

Using the universal property of the $\infty$-category of Pro-objects and the fact that $\zD(\mathsf{Cond}_R)$ is complete, we obtain a functor 
\begin{equation}\label{eqn:pro_to_Dcond}
    \Pro\left(\zD(R)\right) \longrightarrow \mathscr{D}(\cond_R)
\end{equation}
sending a pro-object $X\colon I\to \zD(R)$ to its limit $\lim_i X_i$ computed in condensed objects.\\

We can now state the main result of this section.

\begin{thm}
    \label{thm:Tate_fully_faithful}
    Let $R$ be a ring of finite type over $\mathbb{Z}$.
    The restriction
    \[
        \Tate_{\aleph_0,R}^\mathrm{b} \longrightarrow \mathscr{D}(\cond_R)
    \]
    of the functor \eqref{eqn:pro_to_Dcond} to bounded countable Tate objects  is fully faithful. 
\end{thm}

\begin{proof}
    First observe that the functor \eqref{eqn:pro_to_Dcond} factors through the full subcategory $\zD(\solid_R)$, so that we may as well work in the latter $\infty$-category.
    Let $X$ and $Y$ be in $\Tate_{\ale,R}^\mathrm{b}$.  
    We want to prove that the comparison map 
    \begin{equation}
        \gamma_{X,Y}\colon \Map_{\Tate_{\ale,R}^\mathrm{b}}(X,Y)  \longrightarrow 
        \Map_{\zD(\solid_R)} (\rho(X),\rho(Y))
    \end{equation}
    is an equivalence. By definition of Tate objects, we can pick elementary Tate objects  $V$ and $W$ of which $X$ and $Y$ are respectively retracts. Replacing $V$ and $W$ by appropriate truncations if necessary, we may further assume that both are bounded. Observe that $\gamma_{X,Y}$ is a retract of $\gamma_{V,W}$; since equivalences are stable under retracts, it will suffice to show that $\gamma_{V,W}$ is an equivalence.
    %we may assume without any loss of generality that $V$ and $W$ are elementary Tate objects.
    As an elementary Tate object, $V$ fits in a fiber sequence $V^{\mathrm{p}}\to V\to V^{\mathrm{i}}$ with $V^{\mathrm{p}} \in \Pro_\ale(\Perf_R)$ and $V^{\mathrm{i}}\in \Ind_\ale(\Perf_R)$ by \cite[Corollary 3.4]{hennion:tate}. Since $V$ is bounded, by applying appropriate truncations to the fiber sequence, we may assume without any loss of generality that $V^\mathrm{p}$ and $V^\mathrm{i}$ are also bounded pro and ind-objects. We decompose similarly $W$.
    Consider the induced commutative diagram
    % https://q.uiver.app/#q=WzAsOSxbMCwwLCJcXE1hcChWJyxXJycpIl0sWzEsMCwiXFxNYXAoVicsVykiXSxbMiwwLCJcXE1hcChWJyxXJykiXSxbMiwxLCJcXE1hcChWLFcnJykiXSxbMSwxLCJcXE1hcChWLFcpIl0sWzAsMSwiXFxNYXAoVixXJycpIl0sWzAsMiwiXFxNYXAoVicnLFcnJykiXSxbMiwyLCJcXE1hcChWJycsVycpIl0sWzEsMiwiXFxNYXAoVicnLFcpIl0sWzAsMV0sWzEsMl0sWzIsM10sWzEsNF0sWzAsNV0sWzUsNl0sWzUsNF0sWzQsM10sWzMsN10sWzYsOF0sWzQsOF0sWzgsN11d
    \[\begin{tikzcd}
        {\Map_{\Tate_R}(V^{\mathrm{i}},W^{\mathrm{p}})} & {\Map_{\Tate_R}(V^{\mathrm{i}},W)} & {\Map_{\Tate_R}(V^{\mathrm{i}},W^{\mathrm{i}})} \\
        {\Map_{\Tate_R}(V,W^{\mathrm{p}})} & {\Map_{\Tate_R}(V,W)} & {\Map_{\Tate_R}(V,W^{\mathrm{i}})} \\
        {\Map_{\Tate_R}(V^{\mathrm{p}},W^{\mathrm{p}})} & {\Map_{\Tate_R}(V^{\mathrm{p}},W)} & {\Map_{\Tate_R}(V^{\mathrm{p}},W^{\mathrm{i}})}
        \arrow[from=1-1, to=1-2]
        \arrow[from=1-1, to=2-1]
        \arrow[from=1-2, to=1-3]
        \arrow[from=1-2, to=2-2]
        \arrow[from=1-3, to=2-3]
        \arrow[from=2-1, to=2-2]
        \arrow[from=2-1, to=3-1]
        \arrow[from=2-2, to=2-3]
        \arrow[from=2-2, to=3-2]
        \arrow[from=2-3, to=3-3]
        \arrow[from=3-1, to=3-2]
        \arrow[from=3-2, to=3-3]
    \end{tikzcd}\]
    with exact rows and columns, and the corresponding diagram for mapping spaces in $\zD(\solid_R)$, we see that it suffices to show the result for $V$ and $W$ either in $\Pro^{\mathrm{b}}_\ale(\Perf_R)$ or in $\Ind^{\mathrm{b}}_\ale(\Perf_R)$.
    We treat the case where $V \in \Pro^{\mathrm{b}}_\ale(\Perf_R)$ and $W\in \Ind^{\mathrm{b}}_\ale(\Perf_R)\subseteq \zD^\mathrm{b}(R)$, the other three cases being easier. Using that every object $V_i$ is compact, the proof reduces to the following claim.
\end{proof}

\begin{lem}
    \label{lem:key_computation}
    For $V \colon I\to \Perf_R$ in $\Pro^{\mathrm{b}}_\ale(\Perf_R)$ and $W\in \zD^\mathrm{b}(R)$, the canonical map
    \begin{equation}
        \label{eqn:comparison_pro_to_ind}
        \colim_i \Map_{\zD(\solid_R)}(V_i,W) 
        \longrightarrow 
        \Map_{\zD(\solid_R)}(\lim_i V_i, W)
    \end{equation}
    is an equivalence. 
\end{lem}

\begin{proof}
    Since the cofiltered category $I$ is countable, there exists a coinitial map $\omega^{\mathrm{op}}\to I$ from the opposite category of the ordinal $\omega=(0\to 1\to \dots)$. Without loss of generality, we may therefore suppose that $I\cong \omega^{\mathrm{op}}$. 

    Since $\solid_R$ is an abelian category with enough projectives, we may consider the projective model structure on its category $\mathrm{Ch}^-(\solid_R)$ of right bounded chain complexes. Choose $a$ and $b$ such that each complex $V_i$ is of tor-amplitude in $[a,b-1]$. Using a fibrant replacement,  the pro-object $V$ can be modelled by a diagram $I\to \mathrm{Ch}^-(\solid_R), \, i\mapsto V_i$ in which each $V_i$ is a perfect complex of $R$-modules concentrated in degrees $[a,b]$ and each transition map $V_i \to V_{i-1}$ is an epimorphism.
        
    Let us prove this last claim in details. Starting from any pro-object $i\mapsto V_i$ modelling $V$, we can inductively replace $f_i \colon V_i\to V_{i-1}$ by a map $f'_i \colon V_i' = V_i\oplus X_i \to V_{i-1}$ constructed as follows. Since $V_{i-1}$ is perfect, in every degree $p \in [a,b-1]$ we may choose $n_p\in\mathbb{N}$ and a surjection $g^p\colon R^{n_p} \to V_{i-1}^p$. Consider the standard contractible complex $D(k) = (\dots \to 0\to R\stackrel{\id}{\to} R\to 0 \to \dots)$ with $R$ in (cohomological) degrees $k$ and $k+1$. 
        Observe that there is a unique morphism of complexes $\bar{g}^p\colon D(p)^{\oplus n_p}\to V_{i-1}$ that restricts to $g^p$ in degree $p$. Now take $X_i$ to be $\bigoplus_{p\in [a,b-1]} D(p)^{\oplus n_p}$ and $f'_i$ to be the sum $f\oplus \bigoplus_{p\in [a,b-1]} \bar{g}^p$. By construction, $f_i\colon V_i\to V_{i-1}$ factors through the acyclic cofibration $V_i\stackrel{\sim}{\hookrightarrow} V_i'$, the map $f'_i$ is degreewise surjective and $V_i'$ is a perfect complex concentrated in degrees $[a,b]$, as desired.

    It follows that the diagram $i\mapsto V_i$ we obtain is fibrant in the injective model structure on the category of diagrams $\Fun(I,\mathrm{Ch}^-(\solid_R))$, so that there is a quasi-isomorphism $\lim_i V_i \simeq \holim_i V_i$. We claim that the object $\lim_i V_i$ is actually cofibrant in the projective model structure. Since it is bounded, it suffices to show that it is degreewise projective. Now in degree $p$ the tower is given by surjections $V_i^p\to V_{i-1}^p$ between finitely generated projective modules, which then split, so that the limit $\lim_i V_i^p$ is a direct product $\prod_i Q_i^p$ of finitely generated projective modules $Q_i^p$. Since $R$ is a finitely generated $\mathbb{Z}$-algebra, the product $\prod_\mathbb{N} R$ is projective in $\solid_R$.
        For every $i$ we may write $Q_i^p \oplus P_i^p \cong R^{m_{i,p}}$ for some modules $P_i^p$ and some $m_{i,p}\in \mathbb{N}$, so that we have $\prod_i Q_i^p \oplus \prod_i P_i^p \cong \prod_i R^{m_{i,p}} \cong \prod_\mathbb{N} R$ which implies that $\prod_i Q_i^p$ is also projective. 

        We have shown that $\lim_i V_i$ is cofibrant.
    As a consequence, there is an equivalence of Kan complexes
    \begin{equation}
        \label{eqn:map_as_DK_uhom}
        \Map_{\zD(\solid_R)}(\lim_i V_i, W) \simeq \mathrm{DK}\left(\tau^{\leqslant 0} \uhom_{\mathrm{Ch}^-(\solid_R)}(\lim_i V_i, W) \right)
    \end{equation}
    where $\uhom$ is the $\mathrm{Ch}(R)$-enriched Hom, $\tau^{\leqslant 0} \colon \mathrm{Ch}(R)\to \mathrm{Ch}(R)^{\leqslant 0}$ is the (canonical) truncation functor and
    \[
        \mathrm{DK}\colon \mathrm{Ch}(R)^{\leqslant 0} \simeq \mathsf{sMod}_R \to \mathsf{sSet}
    \]
    is the composite of the Dold--Kan equivalence with the forgetful functor. As every $V_i$ is also cofibrant, formula \eqref{eqn:map_tate_cond} also holds with $\lim_i V_i$ replaced by any of the $V_i$.

    Now in every degree $p$ we have isomorphisms of $R$-modules
    \begin{align*}
        \uhom_{\mathrm{Ch}(\solid_R)}(\lim_i V_i, W)^p 
            & = \bigoplus_{q\in [a,b]} \Hom_{\solid_R}(\lim_i V_i^q, W^{p+q}) \\
            & \cong \bigoplus_{q\in [a,b]} \colim_i \Hom_{\solid_R}(V_i^q, W^{p+q}) \\
            & \cong \colim_i \bigoplus_{q\in [a,b]} \Hom_{\solid_R}(V_i^q, W^{p+q}) \\
            & \cong \colim_i \uhom_{\mathrm{Ch}(\solid_R)}(V_i, W)^p
    \end{align*}
    using the $1$-categorical embedding $\mathsf{Tate}_{\ale,R}\hookrightarrow \solid_R$ (Proposition \ref{ff1}) in the second line. By naturality of these isomorphisms, we obtain an isomorphism of chain complexes 
    \begin{equation*}
        \uhom_{\mathrm{Ch}(\solid_R)}(\lim_i V_i, W)
        \cong 
        \colim_i \uhom_{\mathrm{Ch}(\solid_R)}(V_i, W).
    \end{equation*}
    Finally using that both functors  $\mathrm{DK}$ and $\tau^{\leqslant 0}$ commute with filtered colimits, we deduce that 
    \begin{align*}
        \Map_{\zD(\solid_R)}(\lim_i V_i, W) 
                & \simeq \mathrm{DK}\left(\tau^{\leqslant 0} \uhom_{\mathrm{Ch}(\solid_R)}(\lim_i V_i, W) \right)\\ 
                & \cong \mathrm{DK}\left(\tau^{\leqslant 0} \colim_i \uhom_{\mathrm{Ch}(\solid_R)}(V_i, W) \right)\\ 
                & \cong \colim_i\mathrm{DK}\left(\tau^{\leqslant 0} \uhom_{\mathrm{Ch}(\solid_R)}(V_i, W) \right)\\ 
                & \simeq \colim_i\Map_{\zD(\solid_R)}(V_i, W),
    \end{align*}
    which concludes the proof.
\end{proof}

\begin{rem}
    \label{rem:proof_of_BL}
    A different proof of Lemma \ref{lem:key_computation} can be extracted from Section 2.2 of the recent paper \cite{Bergfalk_Lambie-Hanson}. Although the authors work in the particular case $R=\mathbb{Z}$, one can generalize their arguments to the case where $R$ is of finite type over $\mathbb{Z}$, with only minor changes. Let us sketch this alternative proof for completeness.

    As above, we may assume that $I\cong \omega^{\mathrm{op}}$. The first step consists in writing $\lim V_i$ as an equalizer of two endomorphisms of the product $\prod_i V_i$ to reduce to proving that the morphism of $\infty$-groupoids
    \begin{equation}
        \label{eqn:comporison_colim_to_prod}
        \colim_i   \Map_{\zD(\solid_R)}\left(\bigoplus_{0\leqslant j \leqslant i} V_j, W\right)
        \longrightarrow
        \Map_{\zD(\solid_R)}\left(\prod_i V_i, W\right)
    \end{equation}
    is an equivalence. By dévissage (e.g. using inductively \cite[Proposition 2.13]{Antieau_Gepner_Brauer_group}), since arbitrary products and filtered colimits are exact in $\solid_R$, we may reduce to the case where $V_i$ and $W$ are concentrated in a single degree, with $V_i$ a finitely generated projective $R$-module. Using that $R$ is of finite type over $\mathbb{Z}$ and arguing as in the proof of Lemma \ref{lem:key_computation}, we deduce that $\prod_i V_i$ is projective, therefore both mapping spaces in \eqref{eqn:comporison_colim_to_prod} are discrete. It remains to prove that the map
    \begin{equation*}
        \bigoplus_i  \Hom_{\solid_R}(V_j, W)
        \longrightarrow
        \Hom_{\solid_R}\left(\prod_i V_i, W\right)
    \end{equation*}
    is a bijection; this follows from fully faithfulness of the 1-categorical realization functor $\mathsf{Tate}_{\ale,R} \to \solid_R$ (Proposition \ref{ff1}).
\end{rem}

\begin{rem}
    \label{rem:scholze_counterexample}
    The assumption that $R$ is of finite type over $\mathbb{Z}$ cannot be removed in the statement of Theorem \ref{thm:Tate_fully_faithful}, see \cite{Scholze_mathoverflow} for an explicit counterexample. We thank P. Scholze for pointing it out to us.
% (https://mathoverflow.net/users/6074/peter-scholze)},
\end{rem}

\section{Counter-example to fully faithfulness in the unbounded case}\label{section:counterexample}

In this section, we prove that the boundedness assumption in Theorem \ref{thm:Tate_fully_faithful} cannot be removed.

\begin{prop}
    The functor 
    \begin{equation}
        \label{eqn:functor_Tate_cond}
        \rho\colon \Tate_{\aleph_0}(\Perf_R) \longrightarrow \mathscr{D}(\cond_R)
    \end{equation}
    is not fully faithful.
\end{prop}
\begin{proof}
    We will prove the result by exhibiting countable Tate objects $V$ and $W$ in $\Tate_{\aleph_0}(\Perf_R)$ for which the canonical morphism 
    \begin{equation}
        \label{eqn:map_tate_cond}
        %\gamma \colon 
        \Map_{\Tate_R}(V,W) \longrightarrow 
        \Map_{\zD(\mathsf{Cond}_R)}(\rho (V), \rho(W)).
    \end{equation}
    is not an equivalence. More precisely, we will choose $V$ and $W$ to be respectively in the subcategories $\Pro_{\aleph_0}(\Perf_R)$ and $\Ind_{\aleph_0}(\Perf_R)$ of $\Tate_{\aleph_0}(\Perf_R)$. Writing $V$ as a formal limit   of perfect complexes $V_j$   and using the inclusions $\Ind_{\aleph_0}(\Perf_R)\subseteq \zD(R)\subseteq \zD(\mathsf{Cond}_R)$, we will identify the mapping spaces in \eqref{eqn:map_tate_cond} as 
    \begin{equation*}
        \Map_{\Tate_R}(V,W) \simeq 
        %\colim_j \colim_\ell \Map_{\zD(R)}(V_j, W_\ell) \simeq
        \colim_j \Map_{\zD(R)}(V_j, W)
    \end{equation*}
    and 
    \begin{equation*}
        \Map_{\zD(\mathsf{Cond}_R)}(\rho (V), \rho(W)) \simeq 
        \Map_{\zD(\mathsf{Cond}_R)}(\lim_j V_j, W).
    \end{equation*}
    To prove the result, it will then suffice to show that the comparison map 
    \begin{equation}
        \label{eqn:comparison_map_gamma_bis}
        \gamma\colon  
        \colim_j \Map_{\zD(R)}(V_j, W) \longrightarrow
        \Map_{\zD(\mathsf{Cond}_R)}(\lim_j V_j, W)
    \end{equation}
    is not an equivalence.\\

    We now define the counterexample.
    Let $I$ be the filtered poset $(0 \leftarrow 1 \leftarrow \dots)$ and consider the pro-object $V\colon I\to \zD(\cond_R)$ given by 
    $V_j = \prod_{0 \leq n \leq j} R[n]$, where the transition maps are the obvious projections.
    \begin{claim}
        The limit $\lim_j V_j \simeq \prod_{n\in \mN} R[n]$ is a discrete object in $\zD(\cond_R)$.
    \end{claim}
    \begin{proof}[Proof of the claim]
        Since the inclusion $\zD(R)\to \zD(\cond_R)$ of discrete objects preserves colimits, it suffices to show that the canonical map $\bigoplus_{n\in \mN} R[n]\to \prod_{n\in \mN} R[n]$ is an equivalence. This can be tested at the level of cohomology groups; the result then follows from the fact that $H_n \colon \zD(\cond_R)\to \cond_R$ commutes with direct sums and products for every $n$. %{\color{blue} Verify for product!}
    \end{proof}

    Take $W$ to be the discrete object $\prod_{n\in\mN} R[n]$. %static 
    Observe that the identity morphism $\lim_j V_j \to W$ cannot factor through any of the $V_j$, as any map $V_j\to W$ induces the zero morphism on $H_n$ for $n>j$. This shows that the comparison map $\gamma$ 
    % \begin{equation}
    %     \colim_{j}\Map_{\zD(\cond_R)}(V_j, W) \longrightarrow\Map_{\zD(\cond_R)}\big(\lim_{j} V_j, W\big) 
    % \end{equation}
    described in \eqref{eqn:comparison_map_gamma_bis} is not surjective on connected components, hence not an equivalence.
\end{proof}

\begin{rem}
    Replacing $R[n]$ by $R[-n]$ yields a similar counter-example which instead is left bounded. Therefore neither of the restrictions of the functor \eqref{eqn:functor_Tate_cond} to $\Tate_{\aleph_0}^{[a,\infty)}(\Perf_R)$ nor to $\Tate_{\aleph_0}^{(-\infty,a]}(\Perf_R) $ is fully faithful.
\end{rem}

\begin{rem}
    The counterexample we described above appears independently in \cite[Section 2.2]{Bergfalk_Lambie-Hanson}, where the authors attribute it to L. Mann.
\end{rem}

\bibliographystyle{alpha}
\bibliography{dahema}

\end{document}